\documentclass[12pt]{amsart}

\usepackage{amssymb}
\usepackage{amsfonts}
\usepackage{amsmath}
\usepackage{amssymb}
\usepackage[mathcal]{eucal}
\usepackage{amscd}
\usepackage{verbatim}

\newtheorem{thm}{Theorem}[section]
\newtheorem{cor}[thm]{Corollary}
\newtheorem{lem}[thm]{Lemma}
\newtheorem{prop}[thm]{Proposition}

\theoremstyle{definition}
\newtheorem{defn}[thm]{Definition}

\theoremstyle{remark}

\numberwithin{equation}{section}

\newcommand{\Acal}{\mathcal{A}}

\newcommand{\Dcal}{\mathcal{D}}

\newcommand{\Jcal}{\mathcal{J}}

\newcommand{\Ncal}{\mathcal{N}}

\newcommand{\Xcal}{\mathcal{X}}
\newcommand{\Ycal}{\mathcal{Y}}
\newcommand{\Zcal}{\mathcal{Z}}

\newcommand{\R}{\mathbb{R}}
\newcommand{\C}{\mathbb{C}}
\newcommand{\N}{\mathbb{N}}
\newcommand{\T}{\mathbb{T}}

\newcommand{\Z}{\mathbb{Z}}

\newcommand{\al}{\alpha}
\newcommand{\Ga}{\Gamma}
\newcommand{\ga}{\gamma}
\newcommand{\del}{\delta}

\newcommand{\ep}{\varepsilon}
\newcommand{\sig}{\sigma}

\newcommand{\la}{\lambda}
\newcommand{\La}{\Lambda}
\newcommand{\tet}{\theta}

\newcommand{\br}{\vspace{4 mm}}

\newcommand{\tri}{\bigtriangleup}
\newcommand{\rest}{\upharpoonright}

\newcommand{\ch}{\mathbf{1}}

\newcommand{\cls}{\rm{cls\,}}
\newcommand{\co}{\rm{co\,}}

\newcommand{\Aut}{\rm{Aut\,}}

\newcommand{\Homeo}{\rm{Homeo\,}}

\newcommand{\Id}{\rm{Id}}

\renewcommand{\(}{\bigl(}
\renewcommand{\)}{\bigr)\vphantom{)}}

\newtheorem{rmk}[thm]{Remark}

\newtheorem{exa}[thm]{Example}

\title
{$G$-continuous functions and whirly actions}

\author{E. Glasner and B. Weiss}

\address{Department of Mathematics\\
     Tel Aviv University\\
         Ramat Aviv\\
         Israel}
\email{glasner@math.tau.ac.il}

\address {Institute of Mathematics\\
 Hebrew University of Jerusalem\\
Jerusalem\\
 Israel}
\email{weiss@math.huji.ac.il}

\begin{document}

\begin{abstract}
This paper continues the work \cite{GTW}.
For a Polish group $G$ the notions of $G$-continuous functions and
whirly actions are further exploited to show that:
(i) A $G$-action is whirly iff it admits no
nontrivial spatial factors.
(ii) Every action
of a Polish L\'evy group is whirly.
(iii) There exists a Polish monothetic group which is not
L\'evy but admits a whirly action.
(iv) In the Polish group ${\Aut}(X,\Xcal,\mu)$, for the generic
automorphism $T$ the action of the Polish group
$\La(T)={\cls}\{T^n: n\in \Z\}\subset {\Aut}(X)$ on the
Lebesgue space $(X,\Xcal,\mu)$ is whirly.
(v) The Polish additive group underlying a separable Hilbert space
admits both spatial and whirly faithful actions.
(vi) When $G$ is a non-archimedean Polish group then
every $G$-action is spatial.
\end{abstract}
\maketitle

In the work \cite{GTW} the authors established,
given a boolean action of a Polish group $G$ on
a measure algebra $(\Xcal,\mu)$, a necessary
and sufficient condition for the action $(\Xcal,\mu,G)$
to admit a spatial (or a pointwise) model. This
necessary and sufficient condition was
formulated in terms of certain functions in $L^\infty(\mu)$
called $G$-continuous functions.
Another key notion introduced in \cite{GTW} was that of
a whirly $G$-action. It was shown there that
a whirly action admits no nontrivial spatial factors.
In the present work, which is
a natural continuation of \cite{GTW},
we exploit these new notions and results in several ways.
In the first introductory section we recall the relevant
facts from \cite{GTW}. In the second section we define
stable sets and show that the sub $\sig$-algebra of $\Xcal$
generated by the $G$-continuous functions coincides with the
$\sig$-algebra generated by the stable sets. As corollaries
we deduce that: (i) A $G$-action is whirly iff it admits no
nontrivial spatial factors. (ii) Every action
of a Polish L\'evy group is whirly.
In Section three we deduce some
facts concerning the structure of whirly and spatial systems.
In Section four we show that for a non-archimedean Polish group
every boolean action is spatial.
In the fifth section we show
that in the Polish group ${\Aut}(X,\Xcal,\mu)$ of automorphisms
of a standard Lebesgue space, (topologically) almost every
automorphism defines a whirly action on $(X,\Xcal,\mu)$.
In the final section we show that there are Polish
groups which admit a whirly action and yet
are not L\'evy. First we find $T\in {\Aut}(X)$
for which $\La(T)$ has this property, and then
we show that the abelian Polish group which
underlies a separable infinite dimensional Hilbert space $H$
admits both spatial and whirly faithful actions.
We thank Matt Foreman for a helpful conversation.

\br

\section{Introduction}
In this paper we will
freely use notations and results from \cite{GTW}. The reader
is also referred to that paper for motivation and background.
However for the convenience of the reader we will now
briefly recall some definitions and results from
\cite{GTW}, sometimes in a slightly modified way, that will
be used here.

A \emph{Borel action} of $ G $ on a Borel space
$(X,\Xcal)$ is a Borel map $ G \times X \to X $ satisfying
the conditions, $ ex = x $ and $
g(hx) = (gh)x $ for all $ g,h \in G $ and all $ x \in X $
($e$ denotes the identity element of $G$). Such an
object is called also a \emph{Borel $ G $-space}.
When a $G$-space carries an invariant measure we define
the notion of spatial $G$-action.

\begin{defn}\label{def-spatial}
Let $G$ be a Polish group. By a \emph{spatial $G$-action } we mean a
Borel action of $G$ on a standard Lebesgue space $(X,\Xcal,\mu)$
such that each $g\in G$ preserves the measure $\mu$.  We say that two
spatial actions are isomorphic, if there exists a measure preserving
one to one map between two $G$-invariant subsets of full measure in
the corresponding spaces which intertwines the $G$-actions (the same
two sets for all $g\in G$).
\end{defn}

Often however it is the case that one has to deal with near-actions
(see the definition below and Zimmer \cite[Def.~3.1]{Z})
or merely with an action of the group on a {\em measure
algebra} (i.e. the Borel algebra modulo sets of measure zero) and it is
then desirable to find a spatial model.

\begin{defn}\label{def-near}
Let $G$ be a Polish group and $(X,\Xcal,\mu)$ a standard Borel space
with a probability measure $\mu$.  By a \emph{near-action\/} of $G$ on
$(X,\Xcal,\mu)$ we mean a Borel map $G\times X\to X, (g,x)\mapsto gx$
with the following properties:
\begin{enumerate}
\item[(i)]
With $e$ the identity element of $G$, $ex=x$ for almost every $x$.
\item[(ii)]
For each pair $g,h \in G, g(h x)=(gh) x$ for almost every $x$
(where the set of points $x\in X$ of measure one where this
equality holds may depend on the pair $g, h$).
\item[(iii)]
Each $g\in G$ preserves the measure $\mu$.
\end{enumerate}
\end{defn}

Let $ {\Aut}(X)= {\Aut}(X,\Xcal,\mu) $ be the Polish group of all
equivalence classes of invertible measure preserving transformations
$X \to X$,
with the neighborhood basis at the identity formed by
sets of the form
$$
N(A,\ep)=\{T\in {\Aut}(X): \mu(A\triangle TA)<\ep\},
$$
for $A\in \Xcal$ and $\ep>0$.
The following proposition is from the introduction to
\cite{GTW}.

\begin{prop}
The following three notions
are equivalent.
\begin{enumerate}
\item[(I)]
A near-action of $ G $ on $(X,\Xcal,\mu)$.
\item[(II)]
A continuous homomorphism from $ G $ to $ {\Aut}(X) $.
\item[(III)]
A boolean action of $ G $ on $(X,\Xcal,\mu)$, that is, a continuous
homomorphism from $ G $ to the automorphism group of the associated
measure algebra.
\end{enumerate}
\end{prop}

Every spatial action is also a near-action and in that case the
spatial action will be called a {\em spatial model}
of the near-action.

Recall that a Polish $G$-space is a Polish space $ X $ together
with a continuous action $ G \times X \to X $ of a Polish group $ G
$. Such an action will be called a {\em Polish action}.
If in addition $ X $
is compact then it is a {\em compact} (Polish) $ G $-space.

Every Polish action is also a Borel action. In that case the Polish
action will be called a \emph{Polish model} of the Borel action.

We have the following classical theorems, due to Mackey, Varadarajan
and Ramsay (\cite[Th.~3.2]{Va0}, \cite{Ma},
\cite[Th.~3.3]{Ra} and \cite{Va}).

\begin{thm}\label{lc-model}
Let $G$ be a locally compact second countable topological group.

\textup{(a)}
Every near-action \textup{(}or Boolean action\textup{)} of $ G $
admits a spatial model.

\textup{(b)}
Every spatial action of $ G $ admits a Polish model.
\end{thm}

A powerful generalization to Polish groups of Theorem
\ref{lc-model}(b), given in \cite[Th.~5.2.1]{BK}, is the following.

\begin{thm}[Becker and Kechris]\label{BK-thm}
\mbox{ }

\textup{(a)}
Every Borel action of a Polish group admits a Polish model.

\textup{(b)}
Every Borel $G$-space is embedded \textup{(}as a $G$-invariant Borel
subset\textup{)} into a compact Polish $G$-space.
\end{thm}

It was shown in \cite{GTW} that many Polish groups
have near-actions that do not admit spatial models.
In fact it is shown there that for a Polish L\'evy group
any spatial action is necessarily trivial
(Theorem 1.1 of \cite{GTW}).

L\'evy groups were introduced and studied by
Gromov and Milman in \cite{GM}. We briefly recall the definitions of
L\'evy families and L\'evy groups.
See \cite{MS}, \cite{L}, for detailed
studies of the phenomenon of concentration of measure
and a plethora of examples. Concise proofs of some of the
outstanding instances of this phenomenon can be found in
the appendix of \cite{GTW}.
Refer to \cite{Mi4} and \cite{P1} for further details on
the history of the subject.

A family $(X_n,d_n,\mu_n)$, $n=1,2,3\dots$ of metric spaces
with probability measures $\mu_n$ is called a \emph{L\'evy
family} if the following condition is satisfied.
Whenever $A_n\subset X_n$ is a sequence of subsets such that
$\liminf \mu_n(A_n)>0$ then for any
$\ep >0$, $\lim \mu(B_\ep (A_n))=1$ (here $B_\ep(A)$ is the $\ep$
neighborhood of $A$).
A Polish group $G$ is a \emph{L\'evy group\/} if there exits a family
of compact subgroups $K_n\subset K_{n+1}$ such that the group
$F=\cup_{n\in \N} K_n$ is dense in $G$ and the corresponding family
$(K_n,d,m_n)$ is a L\'evy family; here $m_n$ is the normalized Haar
measure on $K_n$, and $ d $ is a right-invariant compatible metric on
$G$.

Since the group ${\Aut}(X)$ is a L\'evy group it follows
in particular that the natural near-action of ${\Aut}(X)$ on
$(X,\Xcal,\mu)$ does not admit a spatial model.

Another approach used in \cite{GTW} to the problem of
finding spatial models to boolean actions is through the notions of
$G$-continuous function and whirly action.

\begin{defn}
Having a near-action (or boolean action) of $G$ on $(X,\Xcal,\mu)$ we
say that $f \in L^\infty(\mu)$ is \emph{$G$-continuous\/}, if $f \circ
g_n$ converges to $f$ in $L^\infty(\mu)$ norm whenever $g_n \to e$.
\end{defn}

The following theorem is a reformulation of Theorem 2.2
and Remark 2.3 of \cite{GTW}.

\begin{thm}\label{2.2}
Let $(\Xcal,\mu,G)$ be a boolean $G$-system
and let $\Dcal\subset \Xcal$ be the, $G$-invariant,
sub $\sig$-algebra
generated by the $G$-continuous functions
in $L^\infty(\mu)$. Then the factor boolean system
$(\Dcal,\mu,G)$ is the maximum factor of $(\Xcal,\mu,G)$
which admits a spatial model.
Thus $A\in \Xcal$ is in $\Dcal$
iff for every $\ep>0$ there is a $G$-continuous function
$f\in L^\infty$ such that $\|f - \ch_A\|_2< \ep$.
In particular $(\Xcal,\mu,G)$ admits a spatial model
if and only if the collection of $G$-continuous functions
is dense in $L^2(\mu)$.
\end{thm}

The following definition and proposition are from \cite{GTW}
(part (c) of Proposition \ref{3.3} follows directly from
the definition).

\begin{defn}\label{3.2}
A near-action of $G$ on $(X,\Xcal,\mu)$ is \emph{whirly,} if for all
sets $A,B\in\Xcal$ of positive measure, every neighborhood of $e$ (the
unit of $G$) contains $g$ such that $ \mu (A\cap gB) > 0 $.
\end{defn}

\begin{prop}\label{3.3}
\mbox{ }

\textup{(a)}
If a near-action is whirly, then all $G$-continuous functions
are constants.

\textup{(b)}
A whirly action has no spatial model; moreover, such an action cannot
have nontrivial spatial factor.

\textup{(c)}
A factor of a whirly action is whirly.
\end{prop}

In Theorem \ref{cor1} below we will show that conversely, a
boolean action with no nontrivial spatial factor is whirly.

It was shown in \cite{GTW} that the action of ${\Aut}(X)$
on $(X,\Xcal,\mu)$ is whirly, thereby providing a second proof
to the fact that this action does not admit a spatial model.
Theorem \ref{cor2} below asserts that in fact every near-action of a
Polish L\'evy group is whirly.

\br

\section{ Stable subsets }

Let $G$ be a Polish group and fix a sequence
$\{U_m\}_{m=1}^\infty$ of neighborhoods of the
identity element $e\in G$ with the following properties.

(i) For every $m$, $U_{m+1}^2\subset U_m$.

(ii) For every $m$, $U^{-1}_m=U_m$.

(iii) $\bigcap_{m=1}^\infty U_m=\{e\}$.

Let $(\Xcal,\mu,G)$ be a boolean $G$-dynamical system.
We say that $A\in \Xcal$ is {\em positive} when $\mu(A)>0$.
Given a positive set $A\in \Xcal$ and an open set
$U\subset G$, recall that the set $UA \in \Xcal$
was defined in \cite{GTW} as
$$
UA=\bigcup\{\ga A: \ga\in  U'\}.
$$
Here $U'\subset U$ is a countable dense subset of $U$.
This defines $UA$ uniquely as an element of the measure
algebra $\Xcal$ and the definition does not depend on the
choice of $U'$.
Note that for all $u \in U$ we have $uA \subset UA$.

\begin{defn}
For any positive set $A\in \Xcal$ set
$$
\tilde A = \bigcap_{m=1}^\infty U_m A.
$$
We say that $A\in \Xcal$ is {\em stable} if it is
either empty or $\tilde A = A$.
\end{defn}

\begin{prop}\label{finite}
For any $A, B$ positive subsets in $\Xcal$
\begin{enumerate}
\item
$A\subset \tilde A$.
\item
$\tilde{\tilde {A}} = \tilde A.$
\item
For every $g\in G$,
$\widetilde{gA}=g \tilde A$.
\item
$\widetilde{A \cup B}=\tilde A \cup \tilde B.$
\item
$\widetilde{A \cap B}\subset \tilde A \cap \tilde B$, hence for
stable $A,B\in \Xcal$, $\widetilde{A \cap B} = \tilde A \cap
\tilde B$.
\end{enumerate}
\end{prop}

\begin{proof}
1. \
Clear.

2. \
For each $k$,
\begin{align*}
\tilde{\tilde A} & = \bigcap_{m=1}^\infty U_m \tilde A\\
& = \bigcap_{m=1}^\infty U_m
\left(\bigcap_{n=1}^\infty U_n A\right)\\
& \subset U_{k+1}^2 A \subset U_k A,
\end{align*}
hence $\tilde{\tilde A} \subset \bigcap_{k=1}^\infty U_k A
=\tilde A$. By 1, also $\tilde A \subset \tilde{\tilde A}$,
hence $\tilde{\tilde {A}} = \tilde A$.

3.\
Given $g\in G$ and $n\in\N$ there exists an $m\in \N$ with
$g^{-1}U_m g\subset U_n$. Hence
$$
\widetilde{gA}=\bigcap_{m=1}^\infty U_m gA=
\bigcap_{m=1}^\infty g (g^{-1}U_m g) A
\subset g\bigcap_{n=1}^\infty U_n A =g\tilde A.
$$
Now
$$
\tilde A = \widetilde{g^{-1}(g A)} \subset g^{-1}
\widetilde{(g A)},
$$
hence also $g \tilde A\subset \widetilde{gA}$.

The claims 4 and 5 are easy to verify.
\end{proof}

\begin{exa}
Take $G=X=\T$, the unit circle, where $\T$ acts on itself
by translations. Let $A$ be a dense open
subset of measure $1/2$. Then, for $B=X\setminus A$
we have $A\cap B=\emptyset$ but $\tilde A=X $
and $\tilde B =B$, hence $\tilde A \cap \tilde B= B$.
\end{exa}

\begin{prop}\label{intersection}
If $\{A_k\}_{k=1}^\infty$ is a sequence of
stable sets then so is $A=\bigcap_{k=1}^\infty A_k$.
\end{prop}

\begin{proof}
By Proposition \ref{finite} we can assume that
$A_{k+1}\subset A_k$ for every $k$.
Fix $k_0$. Then for each $n$,
$A\subset U_n A \subset U_n A_{k_0}$, hence
$$
\tilde A= \bigcap_{n=1}^\infty U_n A
\subset \bigcap_{n=1}^\infty U_n A_{k_0} =
\tilde A_{k_0}  = A_{k_0}.
$$
Thus
$$
\tilde A \subset
\bigcap_{k=1}^\infty A_k =A.
$$
\end{proof}

\begin{lem}\label{separation}
Let $A$ be a positive stable set and $\ep>0$,
then there exists a stable set $D \subset A^c$
with $\mu(A^c\setminus D)<\ep$ and a $k$
such that $\mu(U_k A \cap U_k D)= 0$.
\end{lem}

\begin{proof}
Let $B=A^c$ and choose an $m$ with $\mu(U_m A\setminus A)<\ep$.
Set $E=(U_mA)^c$, so that $E\subset B$.
If $\mu(U_{m+1}E \cap U_{m+1}A) > 0$, then for some
$\ga, \del \in U_{m+1}$ we have
$\mu(\ga E \cap \del A) > 0$, hence
$\mu(E \cap\ga^{-1} \del A) > 0$.
but $\ga^{-1} \del\in U_{m+1}^2\subset U_m$ and therefore
$\mu(E \cap U_m A) > 0$, contradicting the fact that
$E \subset B= (U_m A)^c$. Thus
$\mu(U_{m+1}E \cap U_{m+1}A) = 0$.

We now set $D=\tilde E$ so that $D$ is stable and
$E \subset  D \subset U_{m+2}E$.
Now $U_{m+2}D \subset U_{m+2}^2E \subset U_{m+1}E$, hence
$\mu(U_k A \cap U_k D)= 0$ with $k=m+2$.
\end{proof}

\begin{prop}\label{complement}
If $A =\tilde A \in \Xcal$ is stable then $A^c$
is a union of an increasing sequence of stable sets.
\end{prop}

\begin{proof}
If $A=\emptyset$ the assertion is clear. Otherwise,
given $\ep$, we use Lemma \ref{separation}
to find a stable subset $D_\ep \subset A^c$ with
$\mu(A^c\setminus D_\ep)< \ep$. Next use an exhaustion
argument to show that $A^c$ is a union of a sequence
of stable sets. Finally we can arrange for this sequence
to be an increasing one by means of Proposition \ref{finite}.
\end{proof}

\begin{defn}
Let $\Acal\subset \Xcal$ be the $\sig$-algebra
generated by the stable sets.
\end{defn}

By Lemma \ref{finite}.3. $\Acal$ is $G$-invariant.

\begin{thm}\label{union}
Every $A \in \Acal$ is a union of an increasing sequence
of stable sets.
\end{thm}

\begin{proof}
1.\
Let $\Acal_0$ be the collection of sets $A\in \Xcal$
which are a union of an increasing sequence
of stable sets. We will show that $\Acal_0 = \Acal$.
Since $\Acal_0$ contains the stable sets it is enough to
show that it is closed under the operations of taking
countable union and taking complements. The closure
under countable union is evident. Suppose
that $A=\bigcup_{k=1}^\infty A_k$ is the union of
the increasing sequence of stable sets $A_k$.
Set $B=A^c$ and $B_k=A_k^c$ so that
$B=\bigcap_{k=1}^\infty B_k$.

By Proposition \ref{complement}
we know that each $B_k$ is an increasing union of stable sets.
Thus, if $\ep>0$ is given we can choose, for each
$k$, a stable subset $D_k\subset B_k$ with
$\mu(B_k\setminus D_k)< \frac{\ep}{2^{k+1}}$.
By Proposition \ref{intersection}
the set  $D_\ep=\bigcap_{k=1}^\infty D_k\subset B$ is stable and
we have
$$
\mu(B\setminus D_\ep)
\le \sum_{k=1}^\infty \mu(B_k\setminus D_k)< \ep.
$$
By an exhaustion argument and an application of
Proposition \ref{finite} we conclude that $B$ is an
increasing union of stable sets.
\end{proof}

The proof of the next theorem is a variation on
the classical Urysohn lemma.

\begin{thm}\label{Ury}
If $(\Xcal,\mu,G)$ is a boolean $G$-system
then $\Dcal=\Acal$; i.e. the sub $\sig$-algebra
generated by the $G$-continuous factors coincides
with the sub $\sig$-algebra generated by the stable sets.
In particular $(\Xcal,\mu,G)$ admits a spatial model
iff $\Xcal=\Acal$.
\end{thm}

\begin{proof}
If $f$ is a $G$-continuous function and
$A_t=\{x: f(x) \le t\}$ then for every $n\in \N$
there exists an $m=m(n)\in \N$ such that
$\|f\circ u - f \|_\infty < 1/n$ for every $u\in U'_m$
(recall that $U'_m$ is a countable dense subset of $U_m$).
For almost every $x\in A_t$ and every $u\in U'_m$ we then have
$$
f(ux) < f(x) +1/n \le t +1/n.
$$
Thus $U_mA_t \subset A_{t+1/n}=\{x: f(x) \le t + 1/n\}$
and it follows that $\tilde A_t = \bigcap_{m=1}^\infty U_mA_t
\subset \bigcap_{n=1}^\infty A_{t+1/n}= A_t$,
so that $\tilde A_t = A_t$. This shows that
$\Dcal \subset \Acal$.

In proving the converse inclusion we can now
restrict our attention to the boolean factor
$(\Acal,\mu,G)$. In other words we now assume,
with no loss of generality, that $\Xcal=\Acal$.

Next observe that, in view of Theorem \ref{2.2},
it suffices to prove the following claim.
Given $A\in \Xcal$ and $\ep>0$ there exists
a $G$-continuous function $f\in L^\infty(\mu)$ such
that $\|f - \ch_A\|_2< \ep$.

Of course we can assume that both $A$ and $A^c$
are positive.

Applying Theorem \ref{union} and
Lemma \ref{separation}
we choose stable subsets $D_0\subset A$, $D_1
\subset A^c$ and an $m_0$ such that
\begin{gather*}
\mu(A\setminus D_0)< \ep/2,\qquad
\mu(A^c\setminus D_1)<\ep/2\\
{\text{and}} \quad \mu(U_{m_0} D_0 \cap U_{m_0} D_1)=0.
\end{gather*}
(We say that $U_{m_0} D_0$ and $U_{m_0} D_1$
are disjoint.)
Define $f$ to have the values $0$ on $D_0$
and $1$ on $D_1$.

Next define, by induction on the binary tree, a system of
disjoint stable subsets
$\{D_r: r=\frac{j}{2^{n}},\ n\in \N,\ 0\le j < 2^n\}$
whose union will be $X$, with the following provision.
Note that if $UD=D$, for a positive set $D$ and a symmetric
neighborhood $U$ of $e$ in $G$, then also $UD^c=D^c$
and both sets are stable.
In the following construction we will assume, for
convenience, that at each stage the various sets
of the form $U_mD_r \setminus D_r$ are positive sets.
If this is not the case then, by the above remark
$D$ as well as $D^c$ are stable and the construction will
terminate at that point. The corresponding interval
of dyadic rationals will be missing from the range of $f$.

\br

Step 1:
Let $E_0=(U_{m_0}(D_0\cup D_1))^c$ and observe that the sets
$U_{m_0+1}E_0$ and $U_{m_0+1}(D_0 \cup D_1)$ are
disjoint. Therefore the stable set
$$
D_{1/2}=\tilde E_0 = \bigl((U_{m_0}(D_0\cup D_1))^c\bigr)^{\sim}
$$
is disjoint from $D_0 \cup D_1$.
Note that
$$
U_{m_0+2}D_{1/2} \subset U_{m_0+2}U_{m_0+2}E_0 \subset
U_{m_0+1}E_0,
$$
hence also the sets $U_{m_0+2}D_0,\ U_{m_0+2}D_{1/2}
 \ {\text{and}} \ U_{m_0+2}D_1$
are pairwise disjoint.
We now choose $m_1 > m_0+2$ such that
\begin{gather*}
\sum\{\mu(U_{m_1}D_r\setminus D_r)
: r=0,1/2,1\}< \frac{1}{2},\\
{\text{and}}\
U_{m_1}D_0,\ U_{m_1}D_{1/2}, \ U_{m_1}D_1\
{\text{are pairwise disjoint}}.
\end{gather*}
We define $f$ to have the value $1/2$ on the stable set $D_{1/2}$
and declare that $f$ will have values in the interval $[0,1/2]$
on $X^1_0=U_{m_0}D_0\cup D_{1/2}$ and values in the interval
$[1/2,1]$ on $X^1_1=U_{m_0}D_1\cup D_{1/2}$, so that
$X^1_0 \cup X^1_1=X$.

\br

Step 2:
Next consider the two disjoint stable sets $D_0$ and
$D_{1/2}$ as subsets of the set $X^1_0=U_{m_0}D_0 \cup D_{1/2}$.
And, at the same time, the
two disjoint stable sets $D_1$ and
$D_{1/2}$ as subsets of the set $X^1_1=U_{m_0}D_1 \cup D_{1/2}$.
Repeating the procedure described in step one we set
$E_1=(U_{m_1}(D_0\cup D_{1/2} \cup D_1))^c$
and note that the sets
$U_{m_1+1}E_1$ and $U_{m_1+1}(D_0 \cup D_{1/2} \cup D_1)$ are
disjoint.
Let
$$
D_{1/4} = \tilde E_1 \cap X^1_0 =
\bigl(U_{m_0}D_0\setminus U_{m_1}(D_0 \cup D_{1/2})\bigr)^{\sim}
$$
and
$$
D_{3/4} = \tilde E_1 \cap X^1_1 =
\bigl(U_{m_0}D_1\setminus U_{m_1}(D_1\cup D_{1/2})\bigr)^{\sim}.
$$
Note that, e.g.,
$$
U_{m_1+2}D_{1/4} \subset U_{m_1+2}U_{m_1+2}E_1 \subset
U_{m_1+1}E_1,
$$
hence also the sets $U_{m_1+2}D_r,\ r\in \{0,1/4,1/2, 3/4,1\}$
are pairwise disjoint.

For a suitable $m_2 > m_1+2$ we will have
\begin{gather*}
\sum\{\mu(U_{m_2}D_r\setminus D_r)
: j=0,1/4,1/2, 3/4,1\}< \frac{1}{2^2},\\
{\text{and the sets}}\
U_{m_2}D_r,\ r \in \{0,1/4,1/2, 3/4,1\}\
{\text{are pairwise disjoint}}.
\end{gather*}
We define $f$ to have the value $r$ on the stable set $D_r$,
$r \in \{0,1/4,1/2, 3/4,1\}$
and declare that $f$ will have values in the interval $[0,1/4]$
on $X^2_{00}=U_{m_1}D_0\cup D_{1/4}$, values in the interval
$[1/4,1/2]$ on $X^2_{01} = (X^1_0\cap U_{m_1}D_{1/2})\cup D_{1/4}$, etc.
Note that we have
$$
U_{m_2}X^1_0 {\text{ is disjoint from $X^2_{11}$
and $U_{m_2}X^1_1$ is disjoint from $X^2_{00}$}}.
$$

\br

Step $n+1$:
Consider, for each $r=r(j)$ in the set
$\{j2^{-(n-1)},\ 0\le j < 2^{(n-1)}\}$
the two disjoint stable sets $D_r$ and $D_{r+2^{-n}}$
as subsets of the set
$X^n_{j0}=(X^{n-1}_{j}\cap U_{m_{n-1}}D_r) \cup D_{r+2^{-n}}$,
and the two disjoint stable sets $D_{r+{2^{-(n-1)}}}$ and $D_{r+2^{-n}}$
as subsets of the set
$X^n_{j1}=(X^{n-1}_{j}\cap U_{m_{n-1}}D_{r+{2^{n-1}}}) \cup D_{r+2^{-n}}$
(here $j0$ denotes the first $n$ bits in the expansion of
$\frac{j}{2^{n-1}}$ in base $2$ and $j1$ the first $n$ bits in
the expansion of $\frac{j}{2^{n-1}} + \frac{1}{2^n}$).
Set
$E_n=\left(\bigcup \{U_{m_n}D_s: s=i2^{-n},\ 0\le i < 2^{n}\}\right)^c$,
then for $r=r(j)=j2^{-(n-1)},\ 0\le j < 2^{(n-1)}$ let
$$
D_{r+ 2^{-(n+1)}} = \tilde E_n \cap X^n_{j0} =
\bigl((X^n_{j0} \cap U_{m_{n-1}}D_r)\setminus
U_{m_n}(D_r \cup D_{r+2^{-n}})\bigr)^{\sim},
$$
and
$$
D_{r+ 2^{-(n-1)} - 2^{-(n+1)}} = \tilde E_n \cap X^n_{j1} =
\bigl((X^n_{j1}\cap U_{m_{n-1}}D_{r+2^{-(n-1)}})\setminus
U_{m_n}(D_{r+2^{-(n-1)}} \cup D_{r+2^{-n}})\bigr)^{\sim}.
$$
For a suitable $m_{n+1} > m_n+2$ we will have
\begin{gather*}
\sum\{\mu(U_{m_{n+1}}D_s\setminus D_s)
: s\in \{i2^{-(n+1)}:\ 0\le i < 2^{(n+1)}\}\}< \frac{1}{2^{n+1}},\\
{\text{and the sets}}\
U_{m_{n+1}}D_s,\  s\in \{i2^{-(n+1)}:\ 0\le i < 2^{(n+1)}\}\
{\text{are pairwise disjoint}}.
\end{gather*}
We define $f$ to have the value $s$ on the stable set $D_s$,
$s\in \{i2^{-(n+1)}:\ 0\le i < 2^{-(n+1)}\}$
and declare that for $r=r(j)=j2^{-(n-1)},\ 0\le j < 2^{(n-1)}$,
$f$ will have values in the interval $[r,r+2^{-(n+1)}]$
on
$$
X^{n+1}_{j00}=(X^{n}_{j0}\cap U_{m_n}D_r)\cup D_{r+2^{-(n+1)}},
$$
and values in the interval $[r+2^{-(n+1)},r+2^{-n}]$ on
$$
X^{n+1}_{j01}=
(X^{n}_{j0}\cap U_{m_n}D_{r+{2^{-n}}})\cup D_{r+2^{-(n+1)}}.
$$
And similarly on the sets
$$
X^{n+1}_{j10}=(X_{j1}\cap U_{m_n}D_{r+2^{-n}})
\cup D_{r+2^{-n}+2^{-(n+1)}}
$$
and
$$
X^{n+1}_{j11}=(X_{j1}\cap U_{m_n}D_{r+2^{-(n-1)}})
\cup D_{r+2^{-n}+2^{-(n+1)}}
$$
$f$ will have values in the intervals $[r+2^{-n},r+2^{-n} +2^{-(n+1)}]$
and $[r+2^{-n} +2^{-(n+1)}, r + 2^{-(n-1)}]$ respectively.

The key point (in proving the $G$-continuity of $f$)
is to observe that, e.g.,
$$
U_{m_{n+1}}X^n_{j0}\subset X^{n+1}_{(j-1)11}
\cup X^n_{j0} \cup X^{n+1}_{j10}.
$$

\br

Once this construction is completed we note that
$$
\mu\left(\bigcup \{D_r:  r=\frac{j}{2^{n}},\ n\in \N,\ 0\le j < 2^n\}
\right)=1,
$$
so that the function
$$
f(x)=\sum_{r} r1_{D_r}
$$
is defined on $X$. It is now easy to conclude that
$\|f - \ch_A\|_2 < \ep$ and that $f$ is $G$-continuous.
\end{proof}

The next two results are direct corollaries of Theorem
\ref{Ury}.

\begin{thm}\label{cor1}
A boolean $G$-system $(\Xcal,\mu,G)$ is whirly
iff it admits no nontrivial spatial factors.
\end{thm}

\begin{proof}
If there are no spatial factors then Theorem \ref{2.2}
implies that there are no nontrivial $G$-continuous functions,
hence, by Theorem \ref{Ury}, no nontrivial stable sets.
It follows that $\mu(UA)=1$ for every positive set
$A\in \Xcal$ and a neighborhood of the identity $U$ in $G$.
By \cite{GTW} this condition is equivalent to $(\Xcal,\mu,G)$
being whirly.

Conversely, if $(\Xcal,\mu,G)$ is whirly then every stable
set is either null or the whole space and by Theorem \ref{Ury},
$\Dcal=\Acal$ is also trivial. Thus $(\Xcal,\mu,G)$
admits no non-constant $G$-continuous functions
and by Theorem \ref{2.2} there are no nontrivial spatial
factors.
\end{proof}

\begin{thm}\label{cor2}
Every ergodic boolean action of a Polish L\'evy group is whirly.
\end{thm}

\begin{proof}
By Theorem 1.1 of \cite{GTW} every ergodic spatial action
of a Polish L\'{e}vy group is trivial; i.e. a one point system.
Evidently our theorem is now a consequence of this fact together
with Theorem \ref{cor1}.
\end{proof}

\br

\section{Some structure theory}
Let $G$ be a Polish group and $\(Y,\Ycal,\nu,G)$ a spatial
$G$-system. Denoting by $M(Y)$ the space  of Borel measures
on the Lebesgue space $Y$, we observe that the spatial action
of $G$ on $Y$ induces an action of $G$ on $M(Y)$. By
\cite{BK} we can choose $(Y,G)$ to be a compact model where
the action $G\times Y \to Y$ is jointly continuous, so that
the induced action $(M(Y),G)$ is jointly continuous as well.
This observation implies that every {\em quasifactor}
(hence in particular every factor) of the system
$\(Y,\Ycal,\nu,G)$ has a spatial model (see \cite{G1}).
Theorem 8.4 of \cite{G1} asserts that two ergodic systems
$(X,\Xcal,\mu,G)$ and $\(Y,\Ycal,\nu,G)$ are {\em not disjoint} iff
$\(Y,\Ycal,\nu,G)$ admits a nontrivial quasifactor which is a factor of
$(X,\Xcal,\mu,G)$.
Thus part 2 of the following theorem is a direct consequence of
Theorem \ref{cor1} above.
For completeness, and for those who are not familiar with
the theory of quasifactors, we provide an alternative direct proof.

\begin{thm}\label{disj}
Let $G$ be a Polish group.
\begin{enumerate}
\item
Every factor of a spatial $G$-system is a spatial system.
\item
A $G$-system is whirly iff it is disjoint from every spatial $G$-system.
\end{enumerate}
\end{thm}

\begin{proof}
1.\
Let $(Z,\Zcal,\eta,G)$ be a boolean factor of the spatial
system $(Y,\Ycal,\nu,G)$.
We consider $\Zcal$ as a $G$-invariant sub $\sig$-algebra
of $\Ycal$ (so that $\eta=\nu\rest \Zcal$)
and let $E^{\Zcal}: L^2(\Ycal,\nu) \to L^2(\Zcal,\nu)$
denote the corresponding conditional expectation; i.e.
the orthogonal projection of the Hilbert space
$L^2(\Ycal,\nu)$ onto the closed subspace $L^2(\Zcal,\nu)$.
Since $E^{\Zcal}$ is a contraction, both on $L^2$ and on $L^\infty$,
it follows easily that if $f\in L^\infty(\Ycal,\nu)$ is a
$G$-continuous function then so is $E^{\Zcal}f$. An application of
Theorem \ref{2.2} yields an $L^2$-dense subset of $L^2(\Ycal,\nu)$
consisting of $G$-continuous functions. Its image under $E^{\Zcal}$
will provide an $L^2$-dense subset of $L^2(\Zcal,\nu)$.
An application of
the other direction of Theorem \ref{2.2} finishes the proof.

\br

2.\
If a $G$-system is disjoint from every spatial
system then in particular it admits no nontrivial spatial
factors and by Theorem \ref{cor1} it is whirly.

Conversely, assume that $(X,\Xcal,\mu,G)$ is
a whirly near-action and that $(Y,\Ycal,\nu,G)$
is spatial.
By \cite{BK} we can assume that $Y$ is a compact space and
that the action $G\times Y \to Y$ is jointly continuous.
Let $\la$ be a joining of the two systems. That is, $\la$ is a
probability measure on $X\times Y$,
with projections $\mu$ and $\nu$ on $X$ and $Y$ respectively.
Let
$$
\la = \int_X (\del_x \times \la_x) \, d\mu(x),
$$ be the disintegration of $\la$ over $\mu$.
Let $\Ga\subset G$ be a countable dense subgroup.
Then, by
the $G$-invariance of $\la$, for
$\mu$-a.e. $x$ we have $\ga \la_x =\ga \la_x$ for all
$\ga\in \Ga$.
If $\la \ne \mu \times \nu$ then the map $x\mapsto \la_x$
is not a constant $\mu$-a.e. and there exists a
continuous function $F$ on $Y$ such that
$$
f(x)=\int_Y F(y)\, d\la_x(y),
$$
is not a constant $\mu$-a.e. Choose a value $c\in \R$
and a $\del >0$
such that
$$
0 < \mu\{x: f(x) > c\} = a < 1,\quad
{\text{and}}\quad
0 < \mu\{x: f(x) > c + \del\} < a.
$$
Let $U$ be a neighborhood of $e$ in $G$ such that
$\sup_{y\in Y}|F(\ga y) - F(y)| < \frac{\del}{10}$
for all $\ga\in U$. Then for $\ga\in U$
$$
f(\ga x) =\int_Y F(y)\, d\la_{\ga x}(y) =
\int_Y F(\ga y)\, d\la_x(y),
$$
for $\mu$-a.e. $x$. Thus
$$
|f(\ga x) - f(x)| \le
\int_Y |F(\ga y) - f(y)|\, d\la_x(y)
\le \frac{\del}{10},
$$
for $\mu$-a.e. $x$.
Denote
$$
A=\{x : f(x) > c+\del\},\qquad B=\{x: f(x) > c\},
$$
then we have for all $\ga\in U$, $\ga A \subset B$, i.e.
$UA \subset B$. This  implies that $A \subset \tilde A \subset B$,
so that $\tilde A$ is a nontrivial stable set in $\Xcal$.
By Theorem \ref{Ury} this contradicts our assumption that
$(X,\Xcal,\mu,G)$ is whirly. Thus $\la=\mu\times \nu$
and we have shown that $(X,\Xcal,\mu,G)$ and $(Y,\Ycal,\mu,G)$
are disjoint.
\end{proof}

\br

Theorem \ref{disj} suggests an analogy ---
whirly is analogous to weak mixing while
spatial is analogous to Kronecker (a system
$(X,\Xcal,\mu,G)$ is a {\em Kronecker system} if
the finite dimensional $G$-subrepresentations are dense
in $L^2(\mu)$).
The following corollary enhances this analogy.

\begin{cor}\label{wm}
Every whirly system is weakly mixing.
\end{cor}
\begin{proof}
Since every Kronecker action embeds in an action of a compact group we
deduce that every Kronecker system admits a spatial model. Since weak
mixing is characterized as the property of having no nontrivial
Kronecker factor, the corollary follows from
Theorem \ref{cor1}
\end{proof}

\br

\section{Non-archimedean group actions are spatial}

Recall that a topological group is called {\em non-archimedean}
if there is a basis for the topology at the identity consisting
of open subgroups. In this section we will apply
our criterion for an action to admit a spatial model
to show that for Polish non-archimedean groups
every near-action admits a spatial model.

Although we will not have an occasion to use it,
we remind the reader of the following interesting
characterization of Polish non-archimedean groups
(see \cite[Theorem 1.5.1]{BK}).

\begin{thm}
A Polish topological group $G$ is non-archimedean iff
it is isomorphic to a closed subgroup of the group $S_\infty$
of permutations of $\N$ (with the topology of pointwise
convergence).
\end{thm}

A theorem that will be used in our proof is the celebrated
Ryll-Nardzewski theorem \cite{RN}
(see also \cite[Theorem III.5.2]{G0}).
Recall that a topological dynamical system $(Q,T)$, where a group $T$
acts continuously on a compact space $Q$, is called {\em affine} if
$Q$ is a convex subset of a topological linear space $E$
and each $t\in T$ acts as an affine transformation
(that is, $t(\al x + (1-\al)y)= \al tx + (1-\al) ty$ for every
$x,y\in Q$ and $0 \le \al \le 1$).
The action is called {\em distal} with respect to a norm $\| \cdot\|$
on $E$ if for every $x\ne y$ in $Q$ we have
$\inf_{t\in T} \|tx - ty\|> 0$.

\begin{thm}\label{rn}
Let $E$ be a separable Banach space,
$Q$ a weakly compact convex subset of $E$. Let $(Q,T)$
be an affine dynamical system such that the action of $T$ is distal
in the norm topology, then $T$ has a fixed point in $Q$.
\end{thm}

\begin{thm}
Every near-action of a Polish non-archimedean group
admits a spatial model.
\end{thm}

\begin{proof}
Let $(X,\Xcal,\mu,G)$ be a near-action of the
Polish non-archimedean group $G$ and let
$g\mapsto U_g$ be the associated Koopman representation
on $L^2(\mu)$ given by $U_g f(x)=f(g^{-1}x),\ g\in G,f\in
L^2(X,\mu)$. Given $f\in L^2(\mu)$ and $\ep>0$ we will next show
that the ball $B_\ep(f)=\{k\in L^2(\mu): \|k - f\|\le \ep\}$
contains a $G$-continuous $L^\infty(\mu)$ function.

By the continuity of the representation $g\mapsto U_g$,
there exists an open subgroup $H\subset G$ such that
$U_h f \in B_\ep(f)$ for every $h\in H$. Being a norm closed
bounded convex subset of $L^2(\mu)$, $B_\ep(f)$ is weakly
closed as well hence weakly compact. We let
$$
Q=w{\text{-}}{\cls}{\co}\left\{ U_h f : h\in H\right\},
$$
the $w$-closed convex hull of the $H$-orbit of $f$ in $L^2(\mu)$.
Clearly $Q$ is a weakly-compact convex $H$-invariant subset
of $B_\ep(f)$. By the Ryll-Nardzewski fixed point
theorem (Theorem \ref{rn})
there exists a function $f_0 \in Q$ which is $H$-fixed.
If we let, for each $M> 0$,
$$
f_M=
  \begin{cases}
   f(x) & \text{if $|f(x)|\le M$}, \\
    -M & \text{if $f(x) < -M$}, \\
    M & \text{if $f(x) > M$},
  \end{cases}
$$
then clearly each $f_M$ is $H$-fixed and
for sufficiently large $M$ we have $f_M\in B_\ep(f)$.
We then fix such an $M$ and claim that $f_M$ is
$G$-continuous.

In fact, the group $G$ being second countable, we see that
the homogeneous space $G/H$ is a countable space
and that consequently the $G$-orbit of $f_M$ in $L^\infty(\mu)$
is a countable set. By Proposition 2.6 of \cite{GTW}
we deduce that $f_M$ is indeed a $G$-continuous function.

Note that what we have shown so far clearly implies
that the $G$-continuous, essentially bounded functions
are dense in $L^2(\mu)$. Thus, in order to complete the proof
of the theorem it only remains to apply Theorem \ref{2.2}.
\end{proof}

\begin{rmk}
One can use a weaker fixed point theorem here. In fact,
since the action of $H$ on $Q$ is {\em weakly almost periodic},
it follows from the general theory of such systems that
every minimal subsystem of $Q$ is equicontinuous
(see e.g. \cite[Chapter 1]{G1}).
Now an affine dynamical system
always contains a {\em strongly proximal} minimal subsystem
(see \cite[Chapter III]{G0}).
Finally, a minimal system which is both proximal and
equicontinuous is necessarily a fixed point.
\end{rmk}

\br

\section{The generic automorphism of a
Lebesgue space is whirly}

Let $(X,\Xcal,\mu)$ be a standard Lebesgue space
and denote by $G={\Aut}(X)$ the Polish group
of its automorphisms, i.e. the equivalence classes
of invertible measure preserving transformations
$X \to X$, equipped with the topology of convergence
in measure. For elements $T$ of $G$ we have the following
well known dichotomy. Either
the subgroup $\{T^n: n\in \Z\}$ is a discrete
subset of $G$ or $\La(T)={\cls}\{T^n: n\in \Z\}$ is a
non-discrete monothetic Polish subgroup of $G$.
In the latter case we say that $T$ is a
{\em rigid} transformation. Thus $T$ is rigid
iff there exists a sequence $n_k\nearrow \infty$
with $\lim_{k\to \infty}T^{n_k}={\Id}$.
A rigid $T$ is called whirly
if the $\La(T)$-action on $(X,\Xcal,\mu)$ is whirly.
This property of $T$ can be characterized directly,
with no reference to $\La(T)$, as follows.

\begin{defn}
An automorphism $T \in G$ is called {\em whirly}
if for every neighborhood $U$ of the identity in $G$
and any two sets $A,B\in \Xcal$ of positive measure
there is $n\in \Z$ such that $T^n\in U$ and
$\mu(T^nA\cap B)>0$.
\end{defn}

The following theorem shows that whirly automorphisms
are not at all rare.

\begin{thm}
The collection of whirly transformations is
a dense $G_\del$ subset of $G$.
\end{thm}

\begin{proof}
For convenience we let $X=[0,1]$ and $\mu$ be normalized
Lebesgue measure on $X$.
Choose a sequence of pairs $\{(A_k, B_k)\}_{k=1}^\infty$
which is dense in $\Xcal \times \Xcal$,
the two-fold product of the measure algebra.
For each $m\in \N$ and $0\le j < 2^m$ let
$J^{(m)}_j=(\frac{j}{2^m},\frac{j+1}{2^m})$ and let
$$
U_m=\{T\in G: \ \mu(T J^{(m)}_j\tri J^{(m)}_j)
< 2^{-2m},\ 0 \le j < 2^m\}.
$$
For $k, m\in \N$ we now define
$$
V_{k,m}=\{T\in G: \exists n \in \Z,\ T^n\in U_m,\
{\text{and}}\ \mu(T^n A_k \cap B_k)
> \del_m\mu(A_k)\mu(B_k)\}.
$$
The sequence of positive constants $\del_m$
will be determined later. Note however that
$\del_m$ is independent of $k$.
Clearly each $V_{k,m}$ is an open subset of $G$
and we set
$$
R=\bigcap_{(k,m)\in \N\times\N}V_{k,m}.
$$

Claim 1: Each $S\in R$ is whirly.

Proof:
Given $A,B$ positive sets in $\Xcal$ and $m_0\in \N$
we have to find an $n\in \N$ with $S^n\in U_m$
and $\mu(S^n A \cap B)>0$. For an $\ep>0$, to be determined soon,
we find a pair $(A_{k_0},B_{k_0})$ with
$$
\mu(A\tri A_{k_0})<\ep, \qquad \mu(B\tri B_{k_0})<\ep.
$$
By assumption $S\in V_{k_0,m_0}$ hence there exists
$n$ such that $S^n \in U_{m_0}$ and
$$
\mu(S^n A_{k_0} \cap B_{k_0})>\del_{m_0}\mu(A_{k_0})
\mu(B_{k_0}).
$$
We now observe that for sufficiently small $\ep$ this
will imply
$$
\mu(S^n A \cap B)>\del_{m_0}\mu(A)\mu(B).
$$

By Baire's category theorem our proof will be complete
once we show that each $V_{k,m}$ is dense in $G$.

Claim 2: For every $m\in \N$ there exists
a number $\del_m> 0$ such that for any positive sets
$A,B \in \Xcal$ the open subset
$$
V=\{T\in G: \exists n \in \Z,\ T^n\in U_m,\
{\text{and}}\ \mu(T^n A \cap B)
> \del_m\mu(A)\mu(B)\}
$$
is dense in $G$.

Proof:
Since the set of weakly mixing automorphisms is
dense in $G$ it suffices to show that, given a weakly
mixing $T_0\in G$ and $\ep>0$, there exists an
element $S\in V$ with $\mu\{x\in X: T_0 x \ne Sx\}<\ep$.

By weak mixing of $T_0$ there exists $n_0$ such that
$\frac{1}{n_0}<< \ep$ and
$\mu(T_0^{n_0} A \cap B)>\frac{1}{2}\mu(A)\mu(B)$.
Denoting $D=T_0^{n_0} A \cap B$ we next apply
the ergodic theorem to choose $n$
sufficiently large so that the frequencies
of appearances of $D$, as well as of each of the
intervals $J^{(m)}_j,\ j=0,\dots,2^m-1$,
in almost every sequence $\{x,T_0x,\dots,T_0^{n-1}x\}$,
is a very good approximation of their respective measures.

Next construct a {\em Rohlin tower} with base $B$,
height $N=nn_0$ and a remainder of measure less than
$\ep$. Thus the sets
$\{B, T_0B, T_0^2 B,\dots, T_0^{N-1}B\}$
are pairwise disjoint and the set $E=X\setminus
\cup_{j=1}^{N-1} T_0^j B$ has measure $< \ep$.

We split $B$ into two parts $B_\ga$
and $B_{1-\ga}$ of measures $\ga\mu(B)$
and $(1-\ga)\mu(B)$ respectively. The size of
$0 < \ga < 1$ will be determined later.

On each of the floors of the $\ga$-part of the tower,
except for the ceiling $T_0^{N-1} B_\ga$,
we set $S=T_0$. On the ceiling we set
$S=T_0^{-N+1}$, so that $S$ becomes periodic
of period $N$ on the $\ga$-part of the tower.

On each of the $n$ $n_0$-block of floors of the
$(\ga-1)$-part of the tower, we define $S$ similarly
to be periodic on that block with period $n_0$,
so that again $S=T_0$ on all floors of the tower except for
the $n$ ceilings of the $n_0$-blocks. Finally, on the remainder
$E$, $S$ is defined as the identity.
Clearly $\mu\{x\in X: T_0 x \ne Sx\}<\ep$ and
we will next show that $S$ is in $V$.

We first observe that if $\ga \le 2^{-2m}$
then $S^{n_0}\in U_m$.
Next note that, by the choice of $n$,
on the $\ga$-part of the tower, say $\tau_\ga$,
we have
$$
\mu(S^{n_0} A\cap B \cap \tau_\ga)\approx \ga\mu(D) >
\frac{\ga}{2}\mu(A)\mu(B),
$$
hence
$$
\mu(S^{n_0} A\cap B) > \frac{\ga}{10}\mu(A)\mu(B).
$$
We are finally in a position to determine the
required size of the constants $\del_m$. If we set
$$
\del_m=\frac{1}{10} 2^{-2m},
$$
then, choosing
$$
\ga= 2^{-2m},
$$
we get
$$
\mu(S^{n_0} A \cap B) >
\frac{\ga}{10}\mu(A)\mu(B) = \frac{2^{-2m}}{10}\mu(A)\mu(B)
= \del_m\mu(A)\mu(B),
$$
so that $S$ is indeed in $V$.
\end{proof}

\br

\section{Further examples}
As was shown in \cite{GTW} a L\'evy group admits no nontrivial
spatial actions. On the other hand every action of
a locally compact group is spatial. In this section
we will show that:
(I) There exists a Polish monothetic group which is not
L\'evy but admits a whirly action.
(II)
The abelian Polish group which
underlies an infinite dimensional
separable Hilbert space $H$ admits both
spatial and whirly faithful actions.

\br

(I)\
In the following discussion we will use some basic
facts from the theory of $IP$-sequences and $IP$-convergence.
We refer to H. Furstenberg's book \cite{F} for the necessary
background.

\begin{defn}
Let $T\in {\Aut}(X)$; we say that the action of $T$ on $X$
is {\em whirly of all orders} if for each $n \ge 1$ the action
of $T\times T \times \cdots\times T$  ($n$-times) on
$X^n$ is whirly.
\end{defn}

\begin{lem}
If for $T\in {\Aut}(X)$ the Polish group $\La(T)$
is L\'evy then the action of $T$ on $X$ is whirly of all
orders.
\end{lem}

\begin{proof}
By Corollary \ref{wm} the system $(X,T)$ is weakly mixing
and therefore for each $n \ge 1$ the system
$(X^n,T\times T \times \cdots\times T)$ is ergodic.
By Theorem \ref{cor2} the latter system is also whirly.
\end{proof}

\begin{lem}
Let $T\in {\Aut}(X)$ be such that the
$T$-action on $X$ is whirly of all orders. Let $U$
be a neighborhood of the identity in $\La(T)$.
Given positive sets $A,B \in \Xcal$ there
exists an $IP$-sequence $\{p_\nu\}$ such that for every $\nu$,
$T^{p_\nu}\in U$ and
$\mu(T^{p_\nu}A \cap B) >0$.
\end{lem}

\begin{proof}
We denote by $\Ncal_e$ the filter of neighborhoods of
the identity in $\La(T)$.
Let $U_1 =U$ and choose $p_1$ such that $T^{p_1}\in U_1$
and $\mu(T^{p_1}A \cap B) >0$.
Let $D=T^{p_1}A \cap B$ and
find $U_2\in \Ncal_e$
such that $U_2\subset U_1$ and $T^{p_1}U_2 \subset U_1$.
Since the action of $T$ is whirly of all orders we can
choose $p_2$ such that $T^{p_2}\in U_2$ and
$\mu(T^{p_2}D \cap B) >0,\
\mu(T^{p_2}A \cap B) >0$. We now have, for $m=p_1, p_2$
and $p_1 + p_2$,
\begin{equation}\label{ip}
T^m \in U \quad {\text{and}}\quad   \mu(T^m A \cap B) >0.
\end{equation}
Assume $p_1,p_2,\dots, p_n$ have been found so that \eqref{ip}
is satisfied for all $m=p_{i_1}+p_{i_2}+\dots + p_{i_k}$
with $i_1 < i_2 < \cdots < i_k \le n$.
For each such $m$ set $D_m = T^m A \cap B$ and choose
$U_{n+1} \in \Ncal_e$ such that $U_{n+1}\subset U_n$ and
$T^m U_{n+1} \subset U$ for every such $m$.
Since the action of $T$ is whirly of all orders
we can choose $p_{n+1}$ with $T^{p_{n+1}}\in U_{n+1}$ and
with the property that
\begin{equation*}
\mu(T^{p_{n+1}} D_m \cap B) >0, \qquad \forall
\ m=p_{i_1}+p_{i_2}+\dots + p_{i_k}.
\end{equation*}
It now follows that \eqref{ip} is satisfied for all
$m=p_{i_1}+p_{i_2}+\dots + p_{i_k}$
with $i_1 < i_2 < \cdots < i_k \le n+1$.
This concludes the inductive step of the construction
of the required $IP$-sequence.
\end{proof}

\begin{prop}\label{char}
Let $T\in {\Aut}(X)$ be such that the
$T$-action on $X$ is whirly of all orders
and let $\al$ be an irrational number. Denote
by $R_\al$ the rotation by $\al$ on $\T$:
$R_\al y = y+\al,  \ y\in \T=\R/\Z$.
Consider the system $(X\times Y, \mu\times m, T\times R_\al)$,
where $m$ is Lebesgue measure on $\T$.
Then
\begin{enumerate}
\item
This dynamical system is rigid and we denote by
$L=\La(T\times R_\al)\subset {\Aut}(X\times \T)$ the corresponding
perfect Polish monothetic topological group.
\item
The action of $L$ on $X$ via the projection on the
first coordinate is whirly.
\item
$L$ admits a continuous character and a nontrivial spatial action.
In particular $L$ is not L\'evy.
\item
Let $K=\{\beta\in \T: ({\Id},\beta)\in L\}$ then
$K$ is a closed subgroup of $\T$ and $K=\{0\}$
iff the natural projection $L \to \La(T)$ is
a topological isomorphism.
\end{enumerate}
\end{prop}

\begin{proof}
1.\
Since $T$ is a rigid transformation there exists a
sequence $n_i\nearrow\infty$ such that $\lim T^{n_i}={\Id}$
in $\La(T)$. It then follows that for some $IP$
sequence $n_\nu$ we have $IP$-$\lim T^{n_\nu}={\Id}$
(see \cite{F}). This fact, in turn, implies that
for some sub-$IP$-sequence $n_{\nu'}$ we have
$IP$-$\lim n_{\nu'} \al=0$ in $\T$ and therefore also
$$
IP{\text{-}}\lim (T^{n_{\nu'}},n_{\nu'} \al)=({\Id},0)
$$
in $L$. Thus $T\times R_\al$ is a rigid transformation
and $L$ is a perfect Polish monothetic group.

2.\
Note first that an arbitrary neighborhood of the identity
in $L$ contains a neighborhood of the form $V=L \cap (U\times
J_\ep)$, where $U$ is a neighborhood of the identity in $\La(T)$
and $J_\ep=\{\beta \in \T: |\beta| < \ep\}$
(\ $|\cdot| $ denotes the distance to the closest integer).
Thus in order to show that the action of $L$ on $X$ is whirly
it suffices to show that for every pair of positive sets
$A,B \in \Xcal$ and every $V$ as above there exists
$g=(S,\beta)\in V$ with $\mu(SA \cap B)>0$.
An application of Lemma \ref{ip} yields an $IP$-sequence
$\{p_\nu\}$ such that for every $\nu$,
$T^{p_\nu}\in U$ and $\mu(T^{p_\nu}A \cap B) >0$.
Since $\T$ is a compact topological group if follows that
for some sub-$IP$-sequence $\{p_{\nu'}\}$,\
$IP{\text{-}}\lim n_{\nu'}\al =0$, hence eventually
$$
(T^{n_{\nu'}},n_{\nu'} \al) \in V = L \cap (U\times J_\ep).
$$
This completes the proof of part 2.

3.\
Clearly the projection $\chi: L \to \T$ to the second coordinate,
$\pi(S,\beta)=\beta$, for $(S,\beta)\in L \subset \La(T) \times \T$
is a continuous homomorphism; i.e. a continuous character.
Via this character $L$ acts spatially on $\T$.
The image of $\chi$ contains the dense subgroup $\{n\al: n\in \Z\}$
and in particular $\chi$ is nontrivial. Since by Theorem 1.1 of
\cite{GTW} a L\'evy group does not admit nontrivial spatial actions
we conclude that $L$ is not L\'evy.

4.\
Straightforward.
\end{proof}

\begin{rmk}
It follows from part 4 of Proposition \ref{char} that
if in addition to the conditions of that proposition
$T$ satisfies also the condition $K=\{0\}$ then
the action of $\La(T)$ on $X$ is whirly yet $\La(T)$ is not
a L\'evy group. Such a $T$ will provide a negative answer
to Problem 2 below (see also \cite[Theorem 3.3]{G}).
\end{rmk}

\br

(II)\
We first define a spatial action of $H$ represented
as $H=\ell_2(\N)$. As usual let $\T=\R/\Z$
be the 1-torus and let $\la$ be normalized
Lebesgue measure on $\T$.
For $x\in \ell_2(\N)$ with coordinates
$x=(x_1,x_2,\dots)$ define a translation
$$
A_1(x): \T^\N\times \T^\N \to \T^\N\times \T^\N,
\qquad
A_1(x)(t,s) = (t+x,s+\sqrt{2}x)  \pmod 1.
$$
Taking $Y=\T^\N\times \T^\N$
and $\mu=\la^\N \times \la^\N$, the map $x\mapsto A_1(x)$
clearly defines a spatial measure preserving system
$(Y,\Ycal,\mu,H)$. Moreover this action is a
topological system on a compact space and
the action is equicontinuous. In fact the image
of $H$ under $A_1: H \to \Homeo(Y)$ is dense in the compact subgroup
of $\Homeo(Y)$
(isomorphic to $\T^\N\times \T^\N$) of translations of $Y$.
Of course the second component is there to make the
action faithful.

\br

We will next describe a faithful whirly near-action of $H$.
This time we choose to view $H$ as the Hilbert space
$L^2([0,1],\la)$.
Start from the Polish monothetic L\'evy group
$G= L_0 \( [0,1], S^1 \) $
of all (equivalence classes of) measurable functions
$ [0,1] \to S^1 $, where $ S^1 = \{ z \in \C : |z| = 1 \} $
and the topology on $G$ is given by $L^2$-convergence
(see \cite{G}).
The map $\tet:H \to G,\ \tet(f)=e^{if}$
defines a continuous homomorphism of $H$ onto $G$.
Consider $G$ as a subgroup of the unitary group $U(H)$
(the elements of $G$ acting as multiplication operators) and
let $(X,\Xcal,\mu,G)$ be the associated Gaussian near-action.
Define for $f\in L^2(\la)$,
$$
A_2(f): X \times X \to X \times X,
\qquad
A_2(f)(x,x') = (\tet(f)x,\tet(\sqrt{2}f)x').
$$
This defines a faithful near-action $(X\times X,\mu\times\mu,H)$.

Let $L={\cls}(A_2(H))$, where the closure is taken in the Polish
group ${\Aut}(X\times X)$. Of course $L\subset G\times G
\subset {\Aut}(X\times X)$, and we claim that actually
$L = G\times G$. In fact, if for each $n$ we denote by
$H_n$ the finite dimensional subspace of $H$ which consists of all
the square integrable functions $f: [0,1] \to \R$ which are
measurable with respect to the partition
$$
\Jcal_n=\left\{[\frac{j}{2^n},\frac{j+1}{2^n}]:
0\le j < 2^n\right\},
$$
and let $G_n = \tet(H_n)$, then it is easily checked that
${\cls} A_2(H_n) = G_n \times G_n$. Since the union of the
increasing sequence of compact groups $G_n\times G_n$
is dense in $G\times G$ we conclude that indeed
$L = G\times G$.

Next observe that the Polish groups $G$ and $G\times G$
are clearly isomorphic (for example via the map
which sends $g\in G$ to the pair $(g_1,g_2)$ of its restrictions to
$[0,1/2]$ and $[1/2,1]$ respectively). In particular,
we deduce that $G\times G$ is a L\'evy group. By
Theorem \ref{cor2} its action on $X\times X$ is whirly.

Now observe that whenever $\La$ is a dense subgroup
of a Polish group $\Ga$ and $(X,\Xcal,\mu,\Ga)$ is a whirly
action then so is the restricted action $(X,\Xcal,\mu,\La)$.
Applying this observation to $A_2(H)\subset {\cls}(A_2(H)) =
G\times G$ we finally conclude that the faithful action
$(X\times X,\mu\times \mu, H)$ via $A_2$ is whirly.

\br

\noindent{\bf Problems:}

1. For a Polish group $G$, is the product of two whirly
$G$-actions whirly?

2. Suppose that the natural action of a closed
subgroup $G\subset {\Aut}(X)$ on $(X,\Xcal,\mu)$
is whirly; is $G$ necessarily a L\'evy group?
In particular, if for $T\in {\Aut}(X)$ the action
of $\La(T)$ on $X$ is whirly, must $\La(T)$ be a
L\'evy group?

3. Is the generic Polish group
$\La(T)={\cls}\{T^n: n\in \Z\}\subset {\Aut}(X)$ L\'evy?

4. For $G = L_0 \( [0,1], S^1 \) $, observe that
the map $A_2: H \to G\times G$, which
is a 1-1 continuous group homomorphism, is not a
homeomorphism. Is there a topologically faithful
whirly action of $H$?

\br

\end{document}